\documentclass[11pt]{amsart}
\usepackage{euscript}
\usepackage[matrix,arrow,curve]{xy}

\newcommand{\Z}{\mathbb{Z}}

\newcommand{\C}{\mathbb{C}}

\newcommand{\iso}{\cong}           
\newcommand{\htp}{\sim}
\newcommand{\CP}[1]{\C {\mathrm P}^{#1}}

\renewcommand{\o}{\omega}

\newtheorem{thm}{Theorem}
\newtheorem{theorem}[thm]{Theorem}

\newtheorem{lemma}[thm]{Lemma}
\newtheorem{prop}[thm]{Proposition}

\newcommand{\acknowledgments}{{\em Acknowledgments.} }

\newcommand{\Fuk}{\EuScript{F}}
\renewcommand{\O}{\EuScript{O}}
\newcommand{\J}{\EuScript{J}}
\newcommand{\Cat}{\EuScript{C}}
\newcommand{\E}{\EuScript{E}}

\begin{document}

\title[Exact Lagrangian submanifolds]{Exact Lagrangian submanifolds in $T^*\!S^n$ \\
  and the graded Kronecker quiver}
\author{Paul Seidel}
\date{January 2004}
\maketitle

\section{Introduction}

The topology of Lagrangian submanifolds gives rise to some of the
most basic, and also most difficult, questions in symplectic
geometry. We have many tools that can be brought to bear on these
questions, but each one is effective only in a special class of
situations, and their interrelation is by no means clear. The present
paper is a piece of shamelessly biased propaganda for a relatively
obscure approach, using Fukaya categories. We will test-ride this
machinery in a particularly simple case, where the computations are
very explicit, and where the outcome can be nicely compared to known
results. For clarity, the level of technical sophistication has been
damped down a little, and therefore the resulting theorem is not
quite the best one can get. Also, in view of the expository nature of
the paper, we do not take the most direct path to the conclusion, but
instead choose a more scenic route bringing the reader past some
classical questions from linear algebra.

Take $M = T^*\!S^n$, $n \geq 2$, with its canonical symplectic
structure. We will consider compact connected Lagrangian submanifolds
$L \subset M$ which satisfy $H^1(L) = 0$ and whose second
Stiefel-Whitney class $w_2(L)$ vanishes. Note that such an $L$ is
automatically orientable (because $w_1(L)$ is the reduction of the
integer-valued Maslov class, which must vanish).

\begin{theorem} \label{th:main}
For such $L$,
\begin{enumerate}
\item \label{item:degree} $[L] = \pm [S^n] \in H_n(M)$;
\item \label{item:cohomology} $H^*(L;\C) \iso H^*(S^n;\C)$;
\item \label{item:fundamental-group}
$\pi_1(L)$ has no nontrivial finite-dimensional complex
representations;
\item \label{item:pair}
if $L' \subset M$ is another Lagrangian submanifold satisfying the
same conditions, then $L \cap L' \neq \emptyset$.
\end{enumerate}
\end{theorem}

As mentioned above, the problem of exact Lagrangian submanifolds in
$M$ and more generally in cotangent bundles has been addressed
previously by several authors, and their results overlap
substantially with Theorem \ref{th:main}. Part \eqref{item:pair},
with somewhat weakened assumptions, was proved by Lalonde-Sikorav
\cite{lalonde-sikorav91} in one of the first papers on the subject,
which is still well worth reading today. For odd $n$, Buhovsky
\cite{buhovsky03} proves a statement essentially corresponding to
\eqref{item:cohomology} (he does not assume $w_2(L) = 0$, and works
with cohomology with $\Z/2$ coefficients; in fact, his result can be
used to remove that assumption from our theorem for odd $n$).
Viterbo's work \cite{viterbo94,viterbo97b} has some relevant
implications, for instance he proved that there is no exact
Lagrangian embedding of a $K(\pi,1)$-manifold in $T^*\!S^n$. Finally,
in the lowest nontrivial dimension $n = 2$ much sharper results are
known. Viterbo's theorem implies that any exact $L \subset T^*S^2$
must be a two-sphere; Eliashberg-Polterovich
\cite{eliashberg-polterovich93} showed that such a sphere is
necessarily differentiably isotopic to the zero-section, and this was
improved to Lagrangian isotopy by Hind \cite{hind03}. Therefore, the
new parts of Theorem \ref{th:main} would seem to be
\eqref{item:degree}, \eqref{item:cohomology} for even $n \geq 4$, and
\eqref{item:fundamental-group} in so far as it goes beyond Viterbo's
results. These overlaps are quite intriguing, as the arguments used
in proving them are widely different (even though all of the
higher-dimensional methods involve Floer homology in some way).

We close with some forward-looking observations. Both Buhovsky's and
Viterbo's approaches should still admit some technical refinements,
as does ours, and one can eventually hope to prove that any exact $L
\subset T^*\!S^n$ is homeomorphic to $S^n$. In contrast, the question
of whether $L$ needs to be diffeomorphic to the standard sphere is
hard to attack, and issues about its Lagrangian isotopy class seem to
be entirely beyond current capabilities (for $n>2$). In the course of
proving Theorem \ref{th:main}, we will show that any submanifold $L$
satisfying its assumptions is ``Floer-theoretically equivalent'' (in
rigorous language, isomorphic in the Donaldson-Fukaya category) to
the zero-section. This is a much weaker equivalence relation than
Lagrangian isotopy, but still sufficiently strong to imply properties
\eqref{item:degree} and \eqref{item:pair} above.

\acknowledgments I would like to thank Paul Biran for fruitful
conversations, and Dan Quillen whose lectures introduced me to the
Kronecker quiver.

\section{The Donaldson-Fukaya category}

In the first place, Fukaya categories provide a convenient way of
packaging the information obtained from Lagrangian Floer cohomology
groups. We will find it convenient to include non-compact Lagrangian
submanifolds with a fixed behaviour at infinity, so we choose a point
on $S^n$ and consider the corresponding cotangent fibre $F_0 \subset
M$. The objects of the Donaldson-Fukaya category $H^0\Fuk(M)$ are
constructed from connected Lagrangian submanifolds $L \subset M$
subject to the following restrictions:
\begin{itemize}
\item
$L$ is either compact, or else agrees with $F_0$ outside a compact
subset.
\item
For any smooth disc with boundary on $L$, $u: (D, \partial D)
\rightarrow (M,L)$, the symplectic area $\int u^*\o$ vanishes;
\item
Similarly, the Maslov number of any disc $u$, which is a relative
Chern class $\int u^*(2c_1(M,L)) \in \Z$, must vanish. As remarked
above, this implies orientability of $L$;
\item $w_2(L) = 0$.
\end{itemize}
To be more precise, the objects are Lagrangian branes $L^\flat$,
which means Lagrangian submanifolds $L$ as described above coming
with certain choices of additional data:
\begin{itemize}
\item
A grading $\tilde{\alpha}_L$, which is a lift of the Lagrangian phase
function $\alpha_L: L \rightarrow S^1$ to a real-valued function.
This exists because of the zero Maslov condition, but clearly there
are infinitely many different choices, differing by integer
constants. We write $L^\flat[k]$ for the brane obtained by shifting
the grading down by $k$.
\item
A spin structure on $L$, and a flat complex vector bundle $\xi_L$.
These are actually coupled, in the sense that changing the spin
structure by an element of $H^1(L;\Z/2)$, and simultaneously
tensoring $\xi_L$ with the corresponding complex line bundle (with
monodromy $\pm 1$), does not change our Lagrangian brane.
\end{itemize}
This list is no doubt baffling to the non-specialist reader, but its
purpose is just to make the Lagrangian Floer cohomology groups
well-defined and as nicely behaved as possibly. We need these to
provide the rest of the category structure, namely, the morphisms
between two objects are the Floer cohomology groups in degree zero,
with twisted coefficients in our flat vector bundles:
\[
 Hom_{H^0\Fuk(M)}(L_0^\flat,L_1^\flat)
 = HF^0(L_0^\flat,L_1^\flat),
\]
and composition of morphisms is provided by the pair-of-pants (or
more appropriately pseudo-holomorphic triangle) product
\begin{equation} \label{eq:triangle-product}
 HF^0(L_1^\flat,L_2^\flat) \otimes HF^0(L_0^\flat,L_1^\flat)
 \longrightarrow HF^0(L_0^\flat,L_2^\flat).
\end{equation}
Some basic reminders about Floer cohomology are in order. First of
all, one can recover the full Floer groups by considering morphisms
into shifted objects, $Hom(L_0^\flat,L_1^\flat[k]) =
HF^k(L_0^\flat,L_1^\flat)$. The endomorphism group of any brane
equipped with the trivial line bundle is its usual cohomology
\begin{equation} \label{eq:floer-is-h}
 HF^*(L^\flat,L^\flat) \iso H^*(L;\C) \quad \text{for} \quad
 \text{$\xi_L \iso \C \times L$.}
\end{equation}
More generally, if we take $L$ with a fixed grading and spin
structure, and equip it with two different flat vector bundles
$\xi_L$, $\xi_L'$, the resulting two branes satisfy
\begin{equation} \label{eq:self-2}
 HF^*(L^{\flat},L^{\flat\,'}) \iso H^*(L;\xi_L^* \otimes \xi_L').
\end{equation}

In our arguments, three simple Lagrangian submanifolds in $M$ will be
prominent. One is the zero-section $Z$, which we make into a brane by
choosing some grading and equipping it with the trivial line bundle.
The second is the fibre $F_0$, treated in the same way, and the third
is the image $F_1 = \tau_Z(F_0)$ under the Dehn twist $\tau_Z$, which
inherits a brane structure from $F_0$.

\begin{lemma} \label{th:floer}
\begin{enumerate}
\item \label{item:floer-1}
The groups $HF^*(Z^\flat,F_k^\flat)$, $HF^*(F_k^\flat,Z^\flat)$ for
$k = 0,1$ are all one-dimensional;
\item \label{item:floer-2}
$HF^*(F_0^\flat,F_1^\flat) \iso H^*(S^{n-1};\C)$;
\item \label{item:floer-3}
$HF^*(F_1^\flat,F_0^\flat) = 0$.
\end{enumerate}
\end{lemma}

Part \eqref{item:floer-1} is obvious, because the Lagrangian
submanifolds intersect in a single point. For the rest, one needs to
remember that Floer cohomology of a pair of non-compact Lagrangian
submanifolds is defined by moving the first one slightly by the
normalized geodesic flow $\phi_t$, which is the Hamiltonian flow of
the function $H$ with $H(\xi) = |\xi|$ outside a compact subset. As a
consequence, the intersections
\[
 \phi_t(L_0) \cap L_1, \quad \text{for $t>0$ small}
\]
are relevant for computing $HF^*(L_0^\flat,L_1^\flat)$. In the case
of $(F_0,F_1)$, the effect of this perturbation is make the
submanifolds intersect cleanly along an $S^{n-1}$, and standard
Bott-Morse methods yield \eqref{item:floer-2}. For $(F_1,F_0)$, in
contrast, the perturbation will make them disjoint, so
\eqref{item:floer-3} follows.

{\em References.} All that we have discussed belongs to the basics of
Floer homology theory. Gradings and twisted coefficients were both
introduced in \cite{kontsevich94}; for the use of spin structures see
\cite{fooo}; for the product \eqref{eq:triangle-product} see
\cite{fukaya-oh98} (these are by no means the only possible
references).

\section{Triangulated categories}

A triangulated category lurks wherever long exact sequences of (any
kind of) cohomology groups appear. The axioms of exact triangles
formalize some nontrivial properties of such sequences, thus allowing
one to manipulate them more efficiently. We will give an informal
description of the axioms, which is neither exhaustive nor totally
rigorous, but which suffices for our purpose. Let $\Cat$ be a
category in which the morphism spaces $Hom(X_0,X_1)$ between any two
objects are complex vector spaces. Having a triangulated structure on
$\Cat$ gives one certain ways of constructing new objects out of old
ones. First of all, we assume that $\Cat$ is additive, which means
that one can form the direct sum $X_0 \oplus X_1$ of two objects,
with the expected properties. Secondly, one can shift or translate an
object by an integer amount, $X \mapsto X[k]$, which allows one to
define higher degree morphism spaces $Hom^k(X_0,X_1) =
Hom(X_0,X_1[k])$. We will denote the direct sum of all these spaces
by $Hom^*(X_0,X_1)$. As a side-remark which will be useful later,
note that by using direct sums and shifts, one can define the tensor
product $V \otimes X$ of any object $X$ with a finite-dimensional
graded complex vector space $V$: choose a homogeneous basis $v_j$ for
$V$, and set
\begin{equation} \label{eq:v-otimes-x}
 V \otimes X = \bigoplus_j V[-deg(v_j)],
\end{equation}
which is easily seen to be independent of the choice of basis up to
isomorphism. The final and most important construction procedure for
objects is the following: to any morphism $a: X_0 \rightarrow X_1$
one can associated a new object, the mapping cone $Cone(a)$, which is
unique up to isomorphism. This can be interpreted as measuring the
failure of $a$ to be an isomorphism: at one extreme, if $a = 0$ then
$Cone(a) = X_0[1] \oplus X_1$, and on the other hand, $a$ is an
isomorphism iff $Cone(a)$ is the zero object.

Mapping cones come with canonical maps $i: X_1 \rightarrow Cone(a)$,
$\pi: Cone(a) \rightarrow X_0[1]$, such that the composition of any
two arrows in the diagram
\begin{equation} \label{eq:triangle}
 X_0 \xrightarrow{a} X_1 \xrightarrow{i} Cone(a) \xrightarrow{\pi}
 X_0[1] \xrightarrow{a[1]} X_1[1]
\end{equation}
is zero. By applying $Hom(Y,-)$ or $Hom(-,Y)$ for some object $Y$ one
gets long exact sequences of vector spaces,
\begin{equation} \label{eq:long-exact}
\begin{aligned}
 & \cdots Hom^k(Y,X_0) \rightarrow Hom^k(Y,X_1) \rightarrow
 Hom^k(Y,Cone(a)) \cdots \\
 & \cdots Hom^k(X_0,Y) \leftarrow Hom^k(X_1,Y) \leftarrow
 Hom^k(Cone(a,Y)) \cdots \\
\end{aligned}
\end{equation}
which is what we were talking about at the beginning of the section.
Diagrams of the form \eqref{eq:triangle} (and isomorphic ones) are
called exact triangles, and drawn rolled up like this:
\[
\xymatrix{ {X_1} \ar[r] & {Cone(a)} \ar[dl]^{[1]} \\
 {X_0} \ar[u]^{a} & }
\]
where the $[1]$ reminds us that this arrow is really a morphism to
$X_0[1]$. One remarkable thing is that triangles can be rotated:
\[
\xymatrix{ {Cone(a)} \ar[r] & {X_0[1]} \ar[dl]^{[1]} \\
 {X_1} \ar[u]^{i} & }
\]
is again an exact triangle, which means that $X_0[1] \iso Cone(i)$,
and similarly $X_1[1] \iso Cone(\pi)$. We know that the cone of the
zero morphism is a direct sum, and therefore in \eqref{eq:triangle}
\begin{equation} \label{eq:splitting}
\begin{aligned}
 a = 0 & \Longrightarrow Cone(a) \iso X_0[1] \oplus X_1, \\
 i = 0 & \Longrightarrow X_0 \iso X_1 \oplus Cone(a)[-1], \\
 \pi = 0 & \Longrightarrow X_1 \iso Cone(a) \oplus X_0. \\
\end{aligned}
\end{equation}
Loosely speaking, the formalism of triangulated categories puts the
boundary operator in long exact sequences on the same footing as the
other two maps. To round off our picture, we look at the situation
where one has two morphisms
\[
 X_0 \xrightarrow{a} X_1 \xrightarrow{b} X_2.
\]
Since the mapping cone measures the failure of a map to be an
isomorphism, it seems intuitive that these defects should somehow add
up under composition, and indeed there is an exact triangle
\begin{equation} \label{eq:octahedral}
\xymatrix{
 {Cone(a)} \ar[r] & {Cone(ba)} \ar[dl]^{[1]} \\
 {Cone(b)[-1].} \ar[u]
}
\end{equation}
This is part of the ``octahedral axiom'' of triangulated structures
(the name comes from a more complicated diagram, which describes
various compatibility conditions between the maps in
\eqref{eq:octahedral} and those in the triangles defining $Cone(a)$,
$Cone(b)$).

The following simple construction arose first in algebraic geometry
as part of the theory of mutations. Suppose that $\Cat$ is a
triangulated category in which the graded vector spaces $V =
Hom^*(X,Y)$ are finite-dimensional. One can define the tensor
products $V \otimes X$ and $V^\vee \otimes Y$ as in
\eqref{eq:v-otimes-x}, and these come with canonical evaluation
morphisms $ev: Hom^*(X,Y) \otimes X \rightarrow Y$, $ev': X
\rightarrow Hom^*(X,Y)^\vee \otimes Y$. We define $T_X(Y) =
Cone(ev)$, $T_Y'(X) = Cone(ev')[-1]$. By definition, this means that
$T_X(Y)$ sits in an exact triangle
\[
\xymatrix{{Y} \ar[r] & {T_X(Y)} \ar[dl]^{[1]} \\
 {Hom^*(X,Y) \otimes X} \ar[u]^{ev} & }
\]
and similarly for $T'_X$. Unfortunately, the axioms of triangulated
categories are not quite strong enough to make $T_X$, $T_X'$ into
actual functors. Still, these operations can be shown to have nice
behaviour, such as taking direct sums to direct sums, and cones to
cones (up to isomorphism).

We will now come to the geometric interpretation of $T,T'$. Contrary
to what the reader may have hoped, the Donaldson-Fukaya category
$H^0\Fuk(M)$ is not triangulated. It does have a shift operation with
the correct properties, namely the change of grading for a Lagrangian
brane, but neither direct sums not cones exist in general. It is an
open question whether one can define a triangulated version of this
category by geometric means, such as including Lagrangian
submanifolds with self-intersections. Meanwhile, there is a purely
algebraic construction in terms of twisted complexes, which yields a
triangulated category $D^b\Fuk(M)$ containing $H^0\Fuk(M)$ as a full
subcategory. The gist of this is simply to add on new objects in a
formal way, so that the requirements of a triangulated category are
satisfied. The details are slightly less straightforward than this
description may suggest, and involve the use of Fukaya's higher order
($A_\infty$) product structures on Floer cohomology, but for the
purposes of this paper, all we need is the knowledge that
$D^b\Fuk(M)$ exists, and the following fact:

\begin{theorem} \label{th:exact}
For any object $L^\flat$ of $H^0\Fuk(M)$, the ``algebraic twist''
$T_{Z^\flat}(L^\flat)$ and the ``geometric (Dehn) twist''
$\tau_Z(L^\flat)$ are isomorphic objects of $D^b\Fuk(M)$. In the same
vein, one has $T_{Z^\flat}'(L^\flat) \iso \tau_Z^{-1}(L^\flat)$. \qed
\end{theorem}

{\em References.} \cite{gelfand-manin} is an accessible presentation
of the abstract theory of triangulated categories. For mutations see
the papers in \cite{rudakov90}. The construction of $D^b\Fuk(M)$ is
outlined in \cite{kontsevich94}. The long exact sequence in Floer
cohomology which is a consequence of Theorem \ref{th:exact} was
introduced in \cite{seidel01}. The argument given there can easily be
adapted to prove the result as stated, see \cite{seidel03b} for an
explanation.

\section{The graded Kronecker quiver}

We will now switch gears slightly. The representation theory of
quivers is a subject with deceptively humble appearance.
Superficially, it is no more than a convenient way of reformulating
certain questions in linear algebra, but in many cases these reveal
themselves to be equivalent to much less elementary problems in other
areas, such as algebraic geometry. We will only need a very simple
instance of the theory, namely the following graded Kronecker quiver,
with nonzero $d \in \Z$:
\begin{equation} \label{eq:kronecker}
\xymatrix{
 {\bullet} \ar@/^0.5pc/[rr]^0 \ar@/_0.5pc/[rr]_d && {\bullet}
}
\end{equation}
By definition, a representation of \eqref{eq:kronecker} consists of
two finite-dimensional graded $\C$-vector spaces $V,W$, together with
linear maps $\alpha: V \rightarrow W$ of degree zero and $\beta: V
\rightarrow W[d]$ of degree $d$, respectively. Two representations
$(V,W,\alpha,\beta)$ and $(V',W',\alpha',\beta')$ are isomorphic if
there are graded linear isomorphisms $\phi: V \rightarrow V'$, $\psi:
W \rightarrow W'$ such that $\alpha'\phi = \psi\alpha$, $\beta'\phi =
\psi\beta$. Of course, any representation splits into a direct sum of
indecomposable ones. The Kronecker quiver is nice in that the latter
can be classified explicitly:

\begin{prop} \label{th:kronecker}
Any indecomposable representation is isomorphic, up to a common shift
in the gradings of both vector spaces, to one of the following:
\begin{itemize}
\item $(\O_X(k),\; k < 0)$: $dim(V) = -k$, $dim(W) = -k-1$,
\[
\xymatrix{
 {V = \C \oplus \C[-d] \oplus \cdots \oplus \C[(k-1)d]}
 \ar@/^0.5pc/[dd]^{\beta = \left(\begin{smallmatrix}
 1 \\ & 1 \\ && \cdots \\ &&& 1 & 0 \end{smallmatrix}\right)}
 \ar@/_0.5pc/[dd]_{\alpha = \left(\begin{smallmatrix}
 0 & 1 \\ && 1 \\ &&& \cdots \\ &&&& 1 \end{smallmatrix}\right)}
 \\
 \\
 {W = \C[-d] \oplus \C[-2d] \oplus \cdots \oplus \C[(k-1)d]}
}
\]
\item $(\O_X(k),\; k \geq 0)$: $dim(V) = k$, $dim(W) = k+1$,
\[
\xymatrix{
 {V = \C[d] \oplus \C[2d] \oplus \cdots \oplus \C[kd]}
 \ar@/^0.5pc/[dd]^{\beta = \left(\begin{smallmatrix}
 1 \\ & 1 \\ && \cdots \\ &&& 1 \\ &&& 0 \end{smallmatrix}\right)}
 \ar@/_0.5pc/[dd]_{\alpha = \left(\begin{smallmatrix}
 0 \\
 1 \\
 & 1 \\
 && \cdots  \\
 &&& 1
 \end{smallmatrix}\right)}
 \\ \\
 {W = \C \oplus \C[d] \oplus \cdots \oplus \C[kd]};
}
\]
\item $(\O_X/\J_{X,0}^k,\;k \geq 1)$: $dim(V) = dim(W) = k$,
\[
\xymatrix{
 {V = \C[d] \oplus \C[2d] \oplus \cdots \oplus \C[kd]}
 \ar@/_0.5pc/[dd]_{\alpha = \left(\begin{smallmatrix}
 0 \\
 1 \\
 & 1 \\
 && \cdots \\
 &&& 1 & 0
 \end{smallmatrix}\right)} \ar@/^0.5pc/[dd]^{\beta = \left(
 \begin{smallmatrix} 1 \\ & 1 \\ && \cdots \\ &&& \cdots \\ &&&&& 1\end{smallmatrix}
 \right)} \\ \\
 {W = \C \oplus \C[d] \oplus \cdots \oplus \C[(k-1)d]}
}
\]
\item $(\O_X/\J_{X,\infty}^k,\;k \geq 1)$: $dim(V) = dim(W) = k$,
\[
\xymatrix{
 {V = \C \oplus \C[d] \oplus \cdots \oplus \C[(k-1)d]}
 \ar@/_0.5pc/[dd]_{\alpha = \left(
 \begin{smallmatrix} 1 \\ & 1 \\ && \cdots \\ &&& \cdots \\ &&&&& 1\end{smallmatrix}
 \right)}
 \ar@/^0.5pc/[dd]^{\beta = \left(\begin{smallmatrix}
 0 & 1 \\
 && 1 \\
 &&& \cdots \\
 &&&& 1 \\
 &&&& 0
 \end{smallmatrix}\right)}
 \\ \\
 {W = \C \oplus \C[d] \oplus \cdots \oplus \C[(k-1)d].}
}
\]
\end{itemize}
\end{prop}

This is actually simple enough to be tackled using only linear
algebra, but the meaning of the classification becomes much clearer
in algebro-geometric terms. To any representation
$(V,W,\alpha,\beta)$ of \eqref{eq:kronecker} one can associate a map
of algebraic vector bundles on $X = \CP{1}$,
\begin{equation} \label{eq:sheaf-map}
 \O_X(-1) \otimes V \xrightarrow{\alpha x + \beta y} \O_X \otimes W.
\end{equation}
Take the $\C^*$-action on $\C^2$ given by $\zeta \cdot (x,y) =
(x,\zeta^{-d}y)$, and the induced action on $X$. Then the sheaves
$\O_X,\O_X(-1)$ are naturally equivariant, and if one equips $V,W$
with the $\C^*$-actions whose weights are given by the grading,
$\alpha x + \beta y$ becomes an equivariant map. At this point, the
usual procedure would be to look at the kernel and cokernel of
\eqref{eq:sheaf-map}, and to use Grothendieck's splitting theorem for
holomorphic vector bundles on $\CP{1}$ to derive Proposition
\ref{th:kronecker}, but we prefer a sleeker and more high-tech
approach using triangulated categories. The category $Coh_{\C^*}(X)$
of equivariant coherent sheaves admits a full embedding into a
triangulated category, its derived category $D^bCoh_{\C^*}(X)$. One
can look at the cone of \eqref{eq:sheaf-map} as an object $\E$ in
that category. A computation using the long exact sequences
\eqref{eq:long-exact} shows that the morphisms between two such cones
are precisely the homomorphisms of quiver representations, which
means that the cone construction embeds the category of
representations of our quiver as a full subcategory into
$D^bCoh_{\C^*}(X)$. In particular, indecomposable representations
must give rise to indecomposable objects $\E$.

Grothendieck's theorem extends easily to the derived category and to
the equivariant case: each indecomposable object of
$D^bCoh_{\C^*}(X)$ is, up to a shift, either a line bundle $\O_X(k)$,
or else a torsion sheaf of the form $\O_X/\J_{X,p}^k$, where $p \in
\{0,\infty\}$ is one of the two fixed points $0 = [0:1]$ or $\infty =
[1:0]$ of the $\C^*$-action. Going back to the objects $\E$
constructed above, one finds that only two essentially different
possibilities can occur. One is that $\E[-1] \iso \O_X(k)$ with
$k<0$, in which case the map \eqref{eq:sheaf-map} is surjective, with
kernel $\E[-1]$. The usual long exact sequence in sheaf cohomology
shows that one can recover the representation from $\E$ as follows:
\[
 V \iso H^1(\E[-1] \otimes \O_X(-1)), \quad
 W \iso H^1(\E[-1]),
\]
and the maps $\alpha,\beta$ are the (Yoneda) products with the
standard generators of $H^0(\O_X(1))$. The other case is where $\E$
is either $\O_X(k)$ with $k \geq 0$, or a torsion sheaf; then
\eqref{eq:sheaf-map} is surjective, with cokernel $\E$, and this time
one finds that
\[
 V \iso H^0(\E \otimes \O_X(-1)), \quad W \iso H^0(\E).
\]
A straightforward computation of cohomology groups identifies the
various $\E$ with the corresponding cases in Proposition
\ref{th:kronecker}, thereby concluding our proof of that result.

{\em References.} The analogue of Proposition \ref{th:kronecker} for
the ungraded quiver is due to Kronecker, and is explained in
textbooks on the representation theory of finite-dimensional algebras
\cite{auslander-reiten-smalo,benson1}. The connection with coherent
sheaves is well-known. The fact that an indecomposable object of the
derived category is actually a single sheaf is a general property of
abelian categories of homological dimension one, and is used
extensively in papers about mirror symmetry on elliptic curves
\cite{polishchuk-zaslow98,kreussler00}. Grothendieck's paper is
\cite{grothendieck57}.

\section{Proof of Theorem \ref{th:main}}

Consider the three basic Lagrangian branes $Z^\flat$, $F_0^\flat$ and
$F_1^\flat$ as objects in $D^b\Fuk(M)$. Because $F_1$ is the image of
$F_0$ under $\tau_Z$, Theorem \ref{th:exact} can be applied, and we
use this to prove:

\begin{lemma} \label{th:kill}
$T_{F_0^\flat}T_{F_1^\flat}(Z^\flat)$ is the zero object.
\end{lemma}

\proof From Lemma \ref{th:floer}\eqref{item:floer-1} we know that
$HF^*(F_1^\flat,Z^\flat)$ is one-dimensional. For simplicity, assume
that the gradings have been chosen in such a way that the nontrivial
Floer group lies in degree zero, and denote a nonzero element of it
by $c$. Using the definitions of $T$, $T'$, and Theorem
\ref{th:exact}, one gets
\[
 T_{F_1^\flat}(Z^\flat) = Cone\big(
 F_1^\flat \xrightarrow{c} Z^\flat\big) = T_{Z^\flat}'(F_1^\flat)[1]
 \iso \tau_Z^{-1}(F_1^\flat)[1] \iso F_0^\flat[1].
\]
By \eqref{eq:floer-is-h} $HF^*(F_0^\flat,F_0^\flat) \iso H^*(F_0;\C)
\iso \C$, and therefore
\[
 T_{F_0^\flat}(F_0^\flat) = Cone\big(
 F_0^\flat \xrightarrow{id} F_0^\flat\big) = 0.
 \qed
\]

Take a Lagrangian submanifold $L \subset M$ satisfying the
assumptions of Theorem \ref{th:main}. Because $H^1(L;\C) = 0$, $L$ is
automatically exact and has zero Maslov class, so we can make it into
an object $L^\flat$ of $H^0\Fuk(M)$ by choosing a grading and spin
structure, as well as the trivial line bundle. Near $Z$, $\tau_Z$ is
equal to the antipodal involution, and so $\tau_Z^2$ is equal to the
identity. By an isotopy along the Liouville (compressing) flow, one
can move $L$ arbitrarily close to $Z$, so $\tau_Z^2(L)$ is Lagrangian
isotopic to $L$. When one considers the gradings on both
submanifolds, however, there is a difference:
\begin{equation} \label{eq:self-shift}
 \tau_Z^2(L^\flat) \htp L^\flat[2-2n].
\end{equation}
From Theorem \ref{th:exact} and the definition of $T$ as a cone, we
get exact triangles
\[
 \xymatrix{
 {L^\flat} \ar[r] & {\tau_Z(L^\flat)} \ar[dl]^{[1]} \\
 {HF^*(Z^\flat,L^\flat) \otimes Z^\flat} \ar[u]
 }
\]
and
\[
 \xymatrix{
 {\tau_Z(L^\flat)} \ar[r] & {\tau_Z^2(L^\flat) \iso L^\flat[2-2n]}
 \ar[dl]^{[1]} \\
 {HF^*(Z^\flat,\tau_Z(L^\flat)) \otimes Z^\flat} \ar[u]
 }
\]
The bottom term can be simplified by noticing that
$HF^*(Z^\flat,\tau_Z(L^\flat)) \iso
HF^*(\tau_Z^{-1}(Z^\flat),L^\flat) \iso HF^*(Z^\flat[n-1],L^\flat)
\iso HF^{*+1-n}(Z^\flat,L^\flat)$. Using the octahedral axiom
\eqref{eq:octahedral}, the two exact triangles can be spliced
together to a single one,
\[
 \xymatrix{
 {L^\flat} \ar[r] & L^\flat[2-2n] \ar[dl]^{[1]} \\
 {Cone\big(HF^*(Z^\flat,L^\flat)[-n] \otimes Z^\flat
 \rightarrow HF^*(Z^\flat,L^\flat) \otimes Z^\flat\big)} \ar[u] }
\]
The $\rightarrow$ must be given by an element of
$HF^{2-2n}(L^\flat,L^\flat)$, but because of \eqref{eq:floer-is-h}
the Floer cohomology groups vanish in negative degrees. Hence the
morphism is necessarily zero, and as explained in
\eqref{eq:splitting},
\begin{equation} \label{eq:split-generated}
 L^\flat \oplus L^\flat[1-2n] \iso
 Cone\big(HF^*(Z^\flat,L^\flat)[-n] \otimes Z^\flat
 \rightarrow HF^*(Z^\flat,L^\flat) \otimes Z^\flat\big).
\end{equation}

\begin{lemma} \label{th:kill-2}
$T_{F_0^\flat}T_{F_1^\flat}(L^\flat)$ is the zero object.
\end{lemma}

\proof From Lemma \ref{th:kill} we know that $T_{F_0^\flat}
T_{F_1^\flat}$ annihilates $Z^\flat$. Since it carries direct sums to
direct sums and cones to cones, this operation also annihilates the
right hand side of \eqref{eq:split-generated}, hence $L^\flat$. \qed

Lemma \ref{th:kill-2} gives us a powerful hold on the a priori
unknown object $L^\flat$. Namely, by putting together the two exact
triangles coming from the definition of $T$ (we omit the details,
since they are parallel to the computation carried out above) one
finds that
\begin{equation} \label{eq:what-is-l}
 L^\flat \iso
 Cone\big( Hom^*(F_0^\flat,T_{F_1^\flat}(L^\flat))[-1] \otimes F_0^\flat
 \rightarrow Hom^*(F_1^\flat, L^\flat) \otimes F_1^\flat
 \big).
\end{equation}
Hence, the isomorphism class of $L^\flat$ as an object of the Fukaya
category is determined by two finite-dimensional graded vector spaces
\[
 V = Hom^*(F_0^\flat,T_{F_1^\flat}(L^\flat))[-1], \quad
 W = Hom^*(F_0^\flat,L^\flat)
\]
and the arrow in \eqref{eq:what-is-l}, which is an element $x \in
Hom_{\C}^*(V,W) \otimes HF^*(F_0^\flat,F_1^\flat)$ of degree zero. We
computed the relevant Floer cohomology group in Lemma
\ref{th:floer}\eqref{item:floer-2}. After choosing generators $a$ of
degree zero and $b$ of degree $n-1$, one writes $x = \alpha \otimes a
+ \beta \otimes b$ where $\alpha: V \rightarrow W$ is a linear map of
degree zero, and $\beta: V \rightarrow W$ one of degree $1-n$. We see
that $(V,W,\alpha,\beta)$ is a representation of the graded Kronecker
quiver with
\begin{equation} \label{eq:d}
d = 1-n < 0.
\end{equation}
If this representation was decomposable, it would give rise to a
corresponding decomposition of the object $L^\flat$ into direct
summands in $D^b\Fuk(M)$, but that is impossible since
$Hom(L^\flat,L^\flat) = HF^0(L^\flat,L^\flat) = H^0(L^\flat;\C) =
\C$. More generally, the long exact sequences \eqref{eq:long-exact}
applied to the cone \eqref{eq:what-is-l} show that $H^*(L;\C) \iso
HF^*(L^\flat,L^\flat)$ is the cohomology of the two-step complex
\begin{multline} \label{eq:end}
 C = \Big\{ 0 \rightarrow Hom_{\C}(V,V) \oplus Hom_{\C}(W,W) \longrightarrow \\
 \xrightarrow{
 \left(\begin{smallmatrix} -\alpha \circ \cdot & \cdot \circ \alpha \\
 -\beta \circ \cdot & \cdot \circ \beta \end{smallmatrix}\right) }
 Hom_{\C}(V,W) \oplus Hom_{\C}(V,W)[d] \rightarrow 0 \Big\}
\end{multline}
This is actually a complex of graded vector spaces, so its cohomology
$H^*(C)$ is bigraded, and one obtains $H^*(L;\C)$ by summing up the
two gradings. One can compute $H^*(C)$ directly for each of the cases
in Proposition \ref{th:kronecker}, or alternatively, one can go back
to our algebro-geometric proof of that result and note that $H^*(C)$
is equal to $Ext^*_X(\E,\E)$ for the associated sheaf $\E$ (which is
bigraded by the cohomological degree and $\C^*$-action). In either
way, one sees that
\begin{itemize}
\item
in the case where the representation is of type $\O_X(k)$, $H^*(C)$
is one-dimensional and concentrated in degree zero;
\item
for $\O_X/\J_{X,0}^k$, $H^*(C)$ is $2k$-dimensional, with generators
in degrees $0,-d,\dots,(1-k)d$ and $d+1,2d+1,\dots,kd+1$.
\item
for $\O_X/\J_{X,\infty}^k$, $H^*(C)$ is $2k$-dimensional, with
generators in degrees $0,d,\dots,(k-1)d$ and $1-d,1-2d,\dots,1-kd$.
\end{itemize}
Bearing in mind \eqref{eq:d}, one sees that $H^*(C)$ cannot be the
cohomology of an $n$-dimensional oriented manifold except in one
case, $\O_X/\J_{X,\infty}$. We have proved:

\begin{lemma} \label{th:final}
Up to a shift, $L^\flat$ is isomorphic to $Cone(a: F_0^\flat
\rightarrow F_1^\flat)$ in the derived Fukaya category. \qed
\end{lemma}

Part \eqref{item:cohomology} of Theorem \ref{th:main} follows
immediately, since in the case of $\O_X/\J_{X,\infty}$, $H^*(C) \iso
H^*(S^n;\C)$. Lemma \ref{th:final} also implies that there is a long
exact sequence
\[
 \cdots HF^*(F_1^\flat,F_1^\flat) \rightarrow
 HF^*(F_1^\flat,L^\flat) \rightarrow HF^{*+1}(F_1^\flat,F_0^\flat) \cdots
\]
which in view of Lemma \ref{th:floer} shows that
$HF^*(F_1^\flat,L^\flat)$ is one-dimensional. Because the Euler
characteristic of Floer cohomology is the ordinary intersection
number, it follows that $F_1 \cdot L = \pm 1$, which implies part
\eqref{item:degree} of Theorem \ref{th:main}. Next, if we had two
Lagrangian submanifolds satisfying the conditions of that theorem,
they would give rise to isomorphic objects in the Fukaya category by
Lemma \ref{th:final}, hence their Floer cohomology would be equal to
the ordinary cohomology of each. Since that is nonzero, the two
Lagrangian submanifolds necessarily intersect, which is part
\eqref{item:pair}.

The remaining statement \eqref{item:fundamental-group} about the
fundamental group is slightly more complicated. Take an
indecomposable flat complex vector bundle $\xi_L$ on $L$, and use
that to define a brane $L^{\flat,'}$. From \eqref{eq:self-2} we see
that $HF^0(L^{\flat,'},L^{\flat,'}) = H^0(\xi_L^* \otimes \xi_L)$
cannot contain any nontrivial idempotents; therefore the
representation of \eqref{eq:kronecker} associated to $L^{\flat,'}$
must still be indecomposable. Again, one finds that
$\O_X/\J_{X,\infty}$ is the only possibility, so $L^{\flat,'}$ is
isomorphic to the previously considered brane $L^\flat$. Now applying
\eqref{eq:self-2} to this pair of branes, one finds that $\xi_L$ must
be the trivial line bundle.

{\em References.} Equation \eqref{eq:self-shift} is taken from
\cite{seidel99}.

\end{document}